 \documentclass{artjlt}            

\usepackage{amssymb}

\title{On the codimension growth of simple color Lie superalgebras}
\author{Du\v san Pagon, Du\v san Repov\v s, and Mikhail Zaicev}
\lastname{Pagon, Repov\v s, and Zaicev}



\address{%
Du\v{s}an~Pagon \\
Faculty of Natural Sciences\\
and Mathematics\\
 University of Maribor\\
P. O. Box 2964\\
Ljubljana, 1001, Slovenia\\
dusan.pagon@uni-mb.si }

\address{%
Du\v san Repov\v s\\
Faculty of Mathematics and Physics\\
University of Ljubljana\\
P.~O.~Box 2964\\
Ljubljana, 1001, Slovenia\\
dusan.repovs@guest.arnes.si
}

\address{%
Mikhail Zaicev \\
Department of Algebra\\
Faculty of Mathematics and Mechanics\\
Moscow State University \\
Moscow,119992, Russia\\
zaicevmv@mail.ru
}

\volume{22}
\annum{2012}

\editor{E. B. Vinberg}
\received{October 5, 2009}

\final{October 3, 2011}

\def\<{\langle}
\def\>{\rangle}
\def\Id{{\rm Id}}

\begin{document}

\maketitle

\begin{abstract}
We study polynomial identities of finite dimensional simple color
Lie superalgebras over an algebraically closed field of
characteristic zero graded by the product of two cyclic groups of
order $2$. We prove that the codimensions  of identities grow
exponentially and the rate of exponent equals  the dimension of
the algebra. A similar result is also obtained for graded identities and
graded codimensions.
\end{abstract}

\section{Introduction}

In this paper we begin to study numerical invariants of polynomial
identities of finite dimensional simple color Lie superalgebras over
an algebraically closed field of characteristic zero. Identities
play an important role in the study of simple algebras. It follows
from the celebrated Amitsur-Levitzky Theorem (see, for example,
\cite[pp.16-18]{GZBook}) that two finite dimensional simple
associative algebras over an algebraically closed field are
isomorphic if and only if they satisfy the same polynomial
identities.  Similar results were later obtained for Lie algebras
\cite{KR}, Jordan algebras \cite{DR} and some other classes. Most
recent results \cite{ShZ} were proved for arbitrary finite
dimensional simple algebras. In the associative case, finite
dimensional graded simple algebras can  also be uniquely defined by
their graded identities \cite{KZ}.

An alternative approach to the characterization of finite dimensional
simple algebras by their identities uses numerical invariants of
identities of algebras. Given an algebra $A$, one can associate with it
 a
sequence of integers $\{c_n(A)\}$, called codimensions of $A$ (all
definitions will be recalled in the next section). If $\dim A=d$,
then it is well-known that $c_n(A)\le d^{n+1}$ (see \cite{GZ1}). For
associative Lie and Jordan algebras it is known that $c_n(A)$ grows
asymptotically like $t^n$, where $t$ is an integer and $0\le t \le
d$ (see  \cite{GRZ}, \cite{GZ1}, \cite{R}). Moreover, $t=d$ if and
only if $A$ is simple.

In the present paper we study the asymptotics of codimensions  of color Lie
superalgebras in the case when $G=\<a\>_2\times \<b\>_2\simeq {\mathbb
Z}_2\oplus {\mathbb Z}_2$ is the product of two cyclic groups of
order two and a skew-symmetric bicharacter $\beta:G\times G
\rightarrow F^*$ is given by $\beta(a,a)=\beta(b,b)=1,
\beta(a,b)=-1$. For any finite dimensional simple Lie algebra $B$,
the corresponding color Lie superalgebra $L=F[G]\otimes B$ is simple
(see \cite{BP}). The main result of the paper asserts
 that the limit
$\lim_{n\to\infty}\sqrt[n]{c_n(L)}$ exists and  equals $\dim A$ (see
Theorem \ref{t1}). All necessary information about polynomial
identities, codimensions and color Lie superalgebras  can be found
in  \cite{BMPZ}, \cite{GZBook}.

\section{Preliminaries}

Let $F$ be a field and  $G$  a finite abelian group. An algebra
$L$ over $F$ is said to be $G$-graded if
$$
L=\bigoplus_{g\in G} L_g
$$
where $L_g$ is a subspace of $L$ and $L_gL_h\subseteq L_{gh}$. An
element $x\in L$ is said to be homogeneous if $x\in L_g$ for some
$g\in G$ and then we say that the degree of $x$ in the grading is
$g$, $\deg x=g$. Any element $x\in L$ can be uniquely decomposed
into a  sum $x=x_{g_1}+\cdots+x_{g_k}$, where $x_{g_1}\in
L_{g_1},\ldots, x_{g_k}\in L_{g_k}$ and $g_1,\ldots,g_k\in G$ are
pairwise distinct. A subspace $V\subseteq L$ is said to be
homogeneous or graded subspace if for any
$x=x_{g_1}+\cdots+x_{g_k}\in V$ we have $x_{g_1},\ldots,x_{g_k}\in
V$. A subalgebra (ideal) $H\subseteq L$ is said to be a graded
subalgebra (ideal) if it is graded as a subspace.

A map $\beta:L\times L\rightarrow F^*$ is said to be a skew-symetric
bicharacter if
$$
\beta(gh,k)=\beta(g,k)\beta(h,k),~\beta(g,hk)=\beta(g,h)\beta(g,k),~
\beta(g,h)\beta(h,g)=1.
$$

A graded $G$-graded algebra $L=\bigoplus_{g\in G} L_g$ is called a
color Lie superalgebra or, more precisely, a $(G,\beta)$-color Lie
superalgebra if for any homogeneous $x,y,z\in L$ one has
$$
xy=-\beta(x,y)yx,
$$
$$
(xy)z=x(yz)-\beta(x,y)y(xz).
$$
Here, for convenience, we write $\beta(x,y)$ instead of $\beta(\deg x,
\deg y)$. Traditionally, the product in color Lie superalgebras is
written as a Lie bracket, $xy=[x,y]$. It is not difficult to see that
$\beta(e,g)=\beta(g,e)=1$, where $e$ is the unit of $G$ and
$\beta(g,g)=\pm 1$ for any $g\in G$ and for any bicharacter $\beta$.
In the case $G=\mathbb Z_2$, $\beta(1,1)=-1$ we get an ordinary Lie
superalgebra. If $\beta(g,g)=1$ for all $g\in G$ then a
$(G,\beta)$-color Lie superalgebra is called  a color Lie algebra.

By definition, a color Lie superalgebra is simple if it has no
non-trivial graded ideals. We study identical relations of
$(G,\beta)$-color Lie algebras in the case $G=\<a\>_2\times
\<b\>_2\simeq {\mathbb Z}_2\oplus {\mathbb Z}_2$  and
$\beta(a,a)=\beta(b,b)=1, \beta(a,b)=-1$. Recently, for these $G$
and $\beta$ all finite dimensional simple color Lie algebras were
classified under a certain weak restriction \cite{BP}. One of the
series of finite dimensional simple algebras can be represented in
the following way.

Let $L=F[G]\otimes B$ be a tensor product of the group ring $F[G]$
with the canonical $G$-grading, and a finite dimensional simple Lie
algebra $B$ with the trivial grading. Then $L$ is a $G$-graded
algebra if we set $\deg(g\otimes x)=g$ for all $g\in G,x\in B$.

Given $i,j,k,l\in \{0,1\}$, we define the product
\begin{equation}\label{1}
[a^ib^j\otimes x, a^kb^l\otimes
y]=(-1)^{j+k}a^{i+k}b^{j+l}\otimes[x,y]
\end{equation}
in $L$. Then under the multiplication (\ref{1}), an algebra $L$ becomes a
$(G,\beta)$-color Lie algebra. Moreover, $L$ is a simple color Lie
algebra.

\begin{remark}\label{r1}
The group algebra $F[G]$ with the multiplication
\begin{equation}\label{1a}
(a^ib^j)\ast(a^kb^l)=(-1)^{j+k}a^{i+k}b^{j+l}
\end{equation}
is isomorphic to $M_2(F)$, the two-by-two matrix algebra over $F$, if we
identify $e,a,b,ab$ with
$$
\left(
   \begin{array}{cc}
    1 & 0  \\
    0 & 1 \\
  \end{array} \right) ,\,
\left(
  \begin{array}{cc}
    -1 & 0  \\
    0 & 1 \\
  \end{array}
  \right),\,
\left(
  \begin{array}{cc}
    0 & 1  \\
    1 & 0 \\
  \end{array}
\right),\, \left(
\begin{array}{cc}
    0 & -1  \\
   1 & 0 \\
  \end{array}
\right),
$$
\end{remark}
\noindent
respectively. \hfill $\Box$ \bigskip

We study non-graded identities of such an algebra $L$.

Next we recall the main notions of the theory of polynomial
identities codimension growth (see \cite{GZBook}). Let $F\{X\}$ be
an absolutely free algebra over $F$ with the countable set of free
generators $X=\{x_1,x_2,\ldots \}$. A non-associative polynomial
$f=f(x_1,\ldots,x_n)$ is said to be an identity of $F$-algebra $A$
if $f(a_1,\ldots,a_n)=0$ for any $a_1,\ldots,a_n\in A$. The set of
all identities of $A$ forms an ideal $\Id(A)$ of $F\{X\}$ stable
under all endomorphisms of $F\{X\}$. Denote by
$P_n=P_n(x_1,\ldots,x_n)$ the subspace of  $F\{X\}$ of all
multilinear polynomials in $x_1,\ldots,x_n$. Then $P_n\cap \Id(A)$ is
a subspace of all multilinear identities of $A$ on variables
$x_1,\ldots,x_n$. A non-negative integer
$$
c_n(A)=\dim\frac{P_n}{P_n\cap \Id(A)}
$$
is called $n$th codimension of $A$. It is well-known
\cite[Proposition 2]{GZ1} that
\begin{equation}\label{ne1}
c_n(A)\le d^{d+1}
\end{equation}
as soon as $\dim A=d<\infty$. In particular, the sequence
$\sqrt[n]{c_n(A)}$ is restricted. In the 1980's, Amitsur conjectured
that the limit $\lim_{n\to\infty}\sqrt[n]{c_n(A)}$ exists and is an
integer for any associative PI-algebra $A$. Amitsur's conjecture was
confirmed for associative \cite{GZ2},\cite{GZ3}, finite dimensional
Lie \cite{Z} and simple special Jordan algebras \cite{GZ1}. For
general non-associative algebras a series of counterexamples with a
fractional rate of exponent were constructed in \cite{GMZ},
\cite{ZM}. If the limit exists we call it the PI-exponent of $A$,
$$
\mbox{PI-exp}(A)=\lim_{n\to\infty}\sqrt[n]{c_n(A)}.
$$

\section{Multialternating polynomials}

Multialternating polynomials play an exceptional role in computing
PI-exponents of simple algebras. In the associative and the Lie case
one may choose multialternating polynomials among central
polynomials constructed by Formanek and Razmyslov. In the Jordan
case the existence of central polynomials is an open problem.
Nevertheless, Razmyslov's approach (see \cite{RAZ}) allows one to
construct the required multialternating polynomials. We shall follow
the Jordan case \cite{GZ4}, \cite{GZ1}.

Recall that $B$ is a finite dimensional simple Lie algebra over an
algebraically closed field of characteristic zero, $G=\<a\>_2\times
\<b\>_2\simeq {\mathbb Z}_2\oplus {\mathbb Z}_2$ and $\beta:G\times
G \rightarrow F^*$ the  skew-symmetric bicharacter on $G$. The
simple color Lie algebra $L$ is equal to $F[G]\otimes B$ and the
multiplication on $L$ is defined by (\ref{1}).

As in the Lie case we define the linear transformation
$\mbox{ad}\,x:L\rightarrow L$ as the right multiplication by $x$,
$\mbox{ad}\,x: y\mapsto [y,x]$. Consider the Killing form $\rho$ on
$L$:
$$
\rho(x,y)= \mbox{tr}(\mbox{ad}\,x\cdot\mbox{ad}\,y).
$$

\begin{lemma}\label{n2}
The Killing form is a symmetric  non-degenerate bilinear form on
$L$.
\end{lemma}
{\em Proof}. Linearity and symmetry of $\rho$ are obvious. Fix any
basis $C=\{c_1,\ldots,c_d\}$ of $B$ where $d=\dim c$ and consider
the basis
$$
\bar C=\{e\otimes c_i, a\otimes c_i, b\otimes c_i, ab\otimes
c_i\vert~1\le i\le d\}
$$
of $L$. Let $M$ be the matrix of $\rho$ in this basis. Consider two
basis elements $x=g\otimes c_i,y=h\otimes c_j\in\bar C$ with $g,
h\in G$. If $g\ne h$ then $gh\ne e$ in $G$ and
$\mbox{ad}\,x\cdot\mbox{ad}\,y$ maps the homogeneous component $L_t$
to $L_{ght}\ne L_t$. Hence
$\mbox{tr}(\mbox{ad}\,x\cdot\mbox{ad}\,y)=0$. Conversely, if $g=h$
then any homogeneous subspace $L_t$ is invariant under the
$\mbox{ad}\,x\cdot\mbox{ad}\,y$-action. Moreover, if we order $\bar
C$ in the following way
$$
\bar C=\{e\otimes c_1,\ldots, e\otimes c_d,a\otimes c_1,\ldots,
a\otimes c_d, b\otimes c_1,\ldots, b\otimes c_d, ab\otimes
c_1,\ldots, ab\otimes c_d
$$
then $M$ is a block-diagonal matrix with four blocks $M_1,\ldots,
M_4$ on the main diagonal and all $M_1,\ldots, M_4$ are matrices of
the Killing form of $B$. Since the Killing form on $B$ is
non-degenerate, the matrix $M$ and $\rho$ are also non-degenerate
and we have completed the proof of the lemma. \hfill $\Box$ \bigskip

Now we fix our simple color Lie algebra $L$, $\dim L=q=4\dim B$ and
construct  multialternating polynomials which are not identities of
$L$. In the rest of this section we shall assume that $F$ is
algebraically closed.

We shall use the following agreement. Given a set of indeterminates
$Y=\{y_1,\ldots,y_n\}$, we denote by $Alt_Y$ the alternation on $Y$.
That is, if \hfill\break $f=f(x_1,\ldots,x_m,y_1,\ldots,y_n)$ is a
polynomial multilinear on $y_1,\ldots,y_n$ then
$$
Alt_Y(f)=\sum_{\sigma \in S_n} (\mbox{sgn}\,\sigma)
f(x_1,\ldots,x_m,y_{\sigma(1)},\ldots,y_{\sigma(n)})
$$
where $S_n$ is the symmetric group and $\mbox{sgn}\,\sigma$ is the
sign of the permutation $\sigma\in S_n$.

\begin{lemma}\label{n0}
Let $B$ be a finite dimensional simple Lie algebra, $\dim L=d$. Then
there exists a left-normed monomial
\begin{equation}\label{b1}
f=[x^1_1,\ldots,x^1_{t_1},y_1,x^2_1,\ldots,x^2_{t_2},y_2,\ldots,
x^d_1,\ldots,x^d_{t_d},y_d,x^{d+1}_1]
\end{equation}
with $t_1,\ldots, t_d>0$ such that $Alt_Y(f)$ is not an identity of
$B$.
\end{lemma}

{\em Proof}. By \cite[Theorem 12.1]{RAZ} there exists a central
polynomial for the pair $(B,\mbox{Ad}\, B)$ that is an associative
polynomial
$$
w=w(x^1_1,\ldots,x^1_{d},\ldots,x^k_1,\ldots,x^k_{d})
$$
such that $w$ is alternating on each set $\{x^i_1,\ldots,x^i_{d}\}$
and
$$
w(\mbox{ad}\,\bar x^1_1,\ldots,\mbox{ad}\,\bar x^k_{d})=\lambda E
$$
is a scalar limear map on $B$ for any evaluation $x^i_j\mapsto\bar
x^i_j\in B$. Moreover, $\lambda\ne 0$ as soon as $\bar
x^i_1,\ldots,\bar x^i_{d}$ are linearly independent for any fixed
$1\le i\le d$.

Hence $[x_0,w]$ is not an identity of $B$. Here we write  $[x_0,w]$
instead of $w(x_0)=w(\mbox{ad}\,\bar x^1_1,\ldots,\mbox{ad}\,\bar
x^k_{d})(x_0)$. By interrupting the alternation on all sets except
$x^k_1,\ldots,x^k_{d}$ and renaming $x^k_1=y_1,\ldots,x^k_{d}=y_d$
we obtain a multilinear polynomial skew-symmetric on
$y_1,\ldots,y_d$ which is not an identity of $B$. By rewriting this
polynomial as a linear combination of left-normed monomials we can
get at least one monomial of the type (\ref{b1}) such that
$Alt_Y(g)$ is not an identity of $B$ but perhaps does not satisfy
the condition $t_1,\ldots,t_d>0$.

If all $t_1,\ldots,t_d>0$ then we are done. Suppose some $t_i=0$.
For brevity we assume $t_1=1,t_2=0$. We again use the central
polynomial. Replace $g=[x^1_1,y_1,y_2,\ldots]$ with
$$
g'=[x^1_1,y_1,w',y_2,\ldots]
$$
where $w'$ is the central polynomial written in new variables
$\widetilde x^i_j$ and we apply $w'(\widetilde x^i_j)$ to
$[x^1_1,y_1]$. Since $w'$ is a central polynomial, one of the
left-normed monomials $f'$ of $g'$ is also of the form (\ref{b1})
with the same $t_1,t_3,\ldots,t_d$ but with $t_2>0$ and $Alt_Y(g')$
is not an identity of $B$. By applying this procedure at most $d$
times we obtain a required polynomial (\ref{b1}). The existence of
the last factor $x^{d+1}_1$ is obvious. \hfill $\Box$\bigskip

Using Lemma \ref{n0} we construct the first alternating polynomial
for $L=F[G]\otimes B$.

\begin{lemma} \label{n3}
There exists a multilinear polynomial
$f=f(x_1,\ldots,x_{q},y_1,\ldots,y_k)$ which is not not vanishing on
$L$ and is alternating on $x_1,\ldots,x_{q}$.
\end{lemma}

{\em Proof}. Let $f$ be the monomial obtained in Lemma \ref{n0}.
Then there exists an evaluation $\varphi: X\rightarrow B$,
$\varphi(x^i_j)=\bar x^i_j, \varphi(y_i)=\bar y_i$, such that
$\varphi(h)\ne 0$ where $h=Alt_Y(f)$. Given $1\le i,j\le 2$, we
consider the evaluation $\varphi _{ij}: X\rightarrow L$ of the
following type:
$$
\varphi_{ij}(y_k)=E_{ij}\otimes \bar y_k,\quad
\varphi_{ij}(x^k_{t_k})=E_{1i}\otimes \bar x^k_{t_k},\quad
\varphi_{ij}(x^{k+1}_1)=E_{j1}\otimes \bar x^{k+1}_{1},\quad 1\le
k\le d,
$$
and
$$
\varphi_{ij}(x^r_{s})=E_{11}\otimes \bar x^r_{s}
$$
for all remaining $x^r_s$ where $E_{ij}$'s are matrix units of
$F[G]\simeq M_2(F)$ (see Remark \ref{r1}). Then
$$
\varphi_{ij}(h)=E_{11}\otimes \varphi(h)\ne 0
$$
in $L$. Now we write $h$ on four disjoint sets of indeterminates,
$$
h_1=h(X_1,Y_1),\ldots,h_4=h(X_4,Y_4).
$$

Since $B$ is simple, the polynomial
$$
H=[h_1,z^1_1,\ldots,z^1_{r_1},h_2,z^2_1,\ldots,z^2_{r_2},\ldots,h_4]
$$
is not an identity of $L$ for some $r_1,\ldots,r_4\ge 0$. Moreover,
\begin{equation}\label{b2}
\varphi_0(Alt(H))= 4d!\cdot [\varphi_{11}(h_1),\bar
z^1_1,\ldots,\bar
z^1_{r_1},\varphi_{12}(h_2),\ldots,\varphi_{22}(h_4)]
\end{equation}
where $\varphi_0\vert_{_{X_1,Y_1}}=\varphi_{11},\ldots,
\varphi_0\vert_{_{X_4,Y_4}}=\varphi_{22}$,
$\varphi_0(z^\gamma_\delta)=\bar z^\gamma_\delta$ and the right hand
side of (\ref{b2}) is non-zero for some $\bar z^\gamma_\delta\in L$.
Here $Alt$ on the left hand side of (\ref{b2}) means the alternation
on $Y_1\cup\ldots\cup Y_4$. Since $|Y_1\cup\ldots\cup Y_4|=4d=\dim
L=q$, we have completed the proof of the lemma. \hfill $\Box$
\bigskip

For extending the number of alternating sets of variables we shall
use the following technical lemma.

\begin{lemma} \label{n4}
Let $f=f(x_1,\ldots, x_m, y_1,\ldots, y_k)$ be a multilinear
polynomial  alternating on $x_1, \ldots, x_m$. Then, for $v,z\in X$,
the polynomial
$$
g=\sum_{i=1}^mf(x_1,\ldots, x_{i-1}, [x_i, v,z], x_{i+1},\ldots,
x_m, y_1,\ldots, y_k)
$$
is also alternating on $x_1, \ldots, x_m$.
\end{lemma}

{\em Proof}.
 Clearly it is enough to check that $g$ is alternating
on $x_r, x_s$, $1\le r<s\le m$. Suppose for instance that $r=1$ and
$s=2$. Since the polynomial
$$
\sum_{i=3}^mf(x_1,\ldots,[x_i, v,z], \ldots, x_m, y_1,\ldots, y_k)
$$
is alternating on $x_1$ and $x_2$, it is enough to check that
$$
g'=f([x_1, v,z],x_2, \ldots, x_m, y_1,\ldots, y_k) + f(x_1, [x_2,
v,z],x_3, \ldots, x_m, y_1,\ldots, y_k)
$$
is alternating on $x_1$ and $x_2$. But
$$
g'(x_1,x_2, \ldots)+g'(x_2,x_1,\ldots) = f([x_1, v,z],x_2, \ldots)
$$
$$
+ f(x_1, [x_2,v,z],\ldots)+ f([x_2,v,z],x_1, \ldots)+ f(x_2, [x_1,
v,z],\ldots)
$$
$$
=f([x_1,v,z],x_2, \ldots)- f([x_2,v,z],x_1, \ldots) +
$$
$$
f([x_2,v,z],x_1, \ldots) - f([x_1,v,z],x_2, \ldots)\equiv 0,
$$
since $f(x,y,\ldots) =-f(y,x,\ldots)$. \hfill $\Box$ \bigskip

In order to simplify the notation, we shall often write \hfill\break
$f=f(x_1,\ldots, x_m, y_1,\ldots, y_n)=f(x_1,\ldots, x_m, Y)$, where
$Y=\{y_1, \ldots, y_n\}$.

\begin{lemma} \label{n5}
Let $Y=Y_0\cup Y_1\cup \cdots \cup Y_r\subseteq X$ be a disjoint
union with $r\ge 0$ and $Y_0$ eventually empty. Let $f=f(x_1,\ldots,
x_{q}, Y)$ be a  multilinear polynomial alternating on each $Y_i, \
1\le i\le r$, and on $x_1, \ldots, x_{q}$. Then, for any $k\ge 1$
and for any $v_1,z_1,\ldots, v_k, z_k\in X$, there exists a
multilinear polynomial
$$
g=g(x_1,\ldots, x_{q}, v_1,z_1,\ldots,  v_k, z_k, Y)
$$
such that, for any evaluation $\varphi: X \to L$, $\varphi(x_i)=
\bar x_i$, $1\le i\le q$, $\varphi(v_j)= \bar v_j$, $\varphi(z_j)=
\bar z_j$, $1\le j\le k$, $\varphi(y)= \bar y$, for $y\in Y$, we
have
$$
\varphi(g)= g(\bar x_1,\ldots, \bar x_{q}, \bar v_1,\bar z_1,\ldots,
\bar v_k, \bar z_k, \bar Y)
$$
$$
= \mbox{tr}(\mbox{ad}\,v_1\cdot\mbox{ad}\,z_1) \cdots
\mbox{tr}(\mbox{ad}\,v_k\cdot \mbox{ad}\,z_k) f(\bar x_1, \ldots,
\bar x_{q}, \bar Y).
$$
Moreover $g$ is alternating on each set $Y_i, \ 1\le i\le r$, and on
$x_1, \ldots, x_{q}$.
\end{lemma}

{\em Proof}. The proof is by induction of $k$. Suppose first that
$k=1$ and define
$$
g=g(x_1,\ldots,  x_{q}, v,z, Y)= \sum_{i=1}^{q} f(x_1, \ldots,
[x_i,v,z], \ldots,x_{q},Y).
$$
Then $g$ is alternating on each set $Y_i, \ 1\le i\le r$ and, by
Lemma \ref{n4}, is also alternating on $x_1, \ldots, x_{q}$.
Consider an evaluation $\varphi: X \to L$ such that $\varphi(x_i)=
\bar x_i$, $1\le i\le q$, $\varphi(v)= \bar v$, $\varphi(z)= \bar
z$, $\varphi(y)= \bar y$, for $y\in Y$. Suppose first that the
elements $\bar x_1,\ldots, \bar x_{q}$ are linearly dependent over
$F$. Then, since $f$ and $g$ are alternating on $x_1,\ldots, x_{q}$,
it follows that $\varphi(f)=\varphi(g)=0$ and we are done.

Therefore we may assume that $\bar x_1,\ldots, \bar x_{q}$ are
linearly independent over $F$ and so,  since $\dim L=q$, they form a
basis of $L$. Hence for all $i=1, \ldots, q$, we write
$$
[\bar x_i \bar v, \bar z] = \alpha_{ii}\bar x_i + \sum_{j\ne i}
\alpha_{ij}\bar x_j,
$$
for some scalars $\alpha_{ij} \in F$. Since $f$ is alternating on
$x_1, \ldots, x_{q}$,
$$
f(\bar x_1, \ldots,[\bar x_i,\bar v, \bar z], \ldots,\bar x_{q},\bar
Y)= \alpha_{ii}f(\bar x_1, \ldots, \bar x_i, \ldots,\bar x_{q},\bar
Y).
$$
Therefore
$$
g(\bar x_1, \ldots,\bar x_{q},\bar v,\bar z,\bar Y) = (\alpha_{11}+
\cdots + \alpha_{qq})f(\bar x_1, \ldots,\bar x_{q},\bar Y),
$$
and, since $\alpha_{11}+ \cdots + \alpha_{qq} =
\mbox{tr}(\mbox{ad}\,v\cdot \mbox{ad}\,z)$, the lemma is thus proved
in case $k=1$.

Now let $k>1$ and let $g=g(x_1,\ldots, x_{q}, v_1,z_1,\ldots,
v_{k-1}, z_{k-1}, Y)$  be a multilinear polynomial satisfying the
conclusion of the lemma. Then we write $g=g(x_1,\ldots, x_{q}, Y')$,
where $Y'=Y_0'\cup Y_1 \cup\cdots\cup Y_r$ and $Y_0'= Y_0\cup\{v_1,
z_1, \ldots, v_{k-1}, z_{k-1}\}$. If we now apply to $g$ the same
arguments as in the case $k=1$, we obtain a polynomial satisfying
the conclusion of the lemma.

\hfill $\Box$\bigskip

Now we are ready to construct the required multialternating
polynomial for our simple color Lie algebra $L$, $\dim L=q$. Recall
that $F$ is an algebraically closed field of characteristic zero.

\begin{proposition} \label{p1}
For any $k\ge 0$ there exists a multilinear polynomial
$$
g_k=g_k(x_1^{(1)},\ldots, x_{q}^{(1)}, \ldots, x_1^{(2k+1)},\ldots,
x_{q}^{(2k+1)}, y_1, \ldots, y_N)
$$
satisfying the following conditions:
\begin{itemize}
\item[1)]
$g_k$ is alternating on each set $\{x_1^{(i)},\ldots,
x_{q}^{(i)}\}$, $1\le i\le 2k+1$;
\item[2)]
$g_k$ is not the identity of $L$;
\item[3)]
the integer $N$ does not depend on $k$.
\end{itemize}
\end{proposition}

{\em Proof}. Let $f=f(x_1,\ldots,x_{q},y_1,\ldots,y_m)$ be the
multilinear polynomial from Lemma \ref{n3}. Hence $f$ is alternating
on $x_1,\ldots,x_{q}$ and does not vanish on $L$.

Suppose first that $k=1$ and write $Y=\{y_1,\ldots, y_m\}$. By Lemma
\ref{n5} there exists a multilinear  polynomial
$$
g=g(x_1,\ldots, x_{q}, v_1^{(1)},z_1^{(1)},\ldots,  v_{q}^{(1)},
z_{q}^{(1)}, Y)
$$
such that
$$
g(\bar x_1, \ldots, \bar x_{q}, \bar v_1^{(1)},\bar
z_1^{(1)},\ldots, \bar v_{q}^{(1)}, \bar z_{q}^{(1)}, \bar Y)
$$
$$
= \mbox{tr}(\mbox{ad}\,\bar v_1^{(1)}\cdot \mbox{ad}\,\bar
z_1^{(1)}) \cdots
 \mbox{tr}(\mbox{ad}\,\bar v_{q}^{(1)} \cdot \mbox{ad}\,\bar z_{q}^{(1)})
 f(\bar x_1, \ldots, \bar x_{q}, \bar Y).
$$
Now, for any $\sigma,\tau\in S_{q}$, define the polynomial
$$
g_{\sigma,\tau}=g_{\sigma,\tau}(x_1, \ldots, x_{q},
v_1^{(1)},z_1^{(1)},\ldots,  v_{q}^{(1)}, z_{q}^{(1)}, Y)
$$
$$
= g(x_1, \ldots, x_{q},
v_{\sigma(1)}^{(1)},z_{\tau(1)}^{(1)},\ldots, v_{\sigma(q)}^{(1)},
z_{\tau(q)}^{(1)}, Y).
$$
Then set
$$
g_1(x_1,\ldots, x_{q}, v_1^{(1)},z_1^{(1)},\ldots, v_{q}^{(1)},
z_{q}^{(1)}, Y)= \frac{1}{q!}\sum_{\sigma, \tau\in S_{q}}
(\mbox{sgn}\, \sigma) (\mbox{sgn}\, \tau) g_{\sigma,\tau}.
$$
The polynomial $g_1$ is alternating on each of the sets $\{x_1,
\ldots x_{q}\}$, $\{v_1^{(1)}, \ldots, v_{q}^{(1)}\}$ and
$\{z_1^{(1)}, \ldots, z_{q}^{(1)}\}$.  Next we show that for any
evaluation $\varphi$, $$\varphi(g_1)= \mbox{det}\bar\rho_1\cdot
\varphi(f),$$ where
$$
\bar \rho_1= \left(
  \begin{array}{ccc}
    \rho(\bar v_1^{(1)}, \bar z_1^{(1)}) & \cdots & \rho(\bar v_1^{(1)}, \bar z_{q}^{(1)}) \\
    \vdots &  & \vdots \\
    \rho(\bar v_{q}^{(1)}, \bar z_1^{(1)}) & \cdots & \rho(\bar v_{q}^{(1)}, \bar z_{q}^{(1)}) \\
  \end{array}
\right).
$$
By Lemma \ref{n5},
$$
\varphi(g_1)= \gamma \varphi(f)
$$
for any evaluation $\varphi: X \to L$, where
$$
\gamma = \frac{1}{q!}\sum_{\sigma, \tau\in S_{q}} (\mbox{sgn}\,
\sigma) (\mbox{sgn}\, \tau) \rho(\bar v_{\sigma(1)}^{(1)}, \bar
z_{\tau(1)}^{(1)}) \cdots \rho(\bar v_{\sigma(q)}^{(1)}, \bar
z_{\tau(q)}^{(1)}).
$$
We fix $\sigma\in S_{q}$ and we compute the sum
$$
\gamma_\sigma = \sum_{\tau \in S_{q}}(\mbox{sgn}\, \tau) \rho(\bar
v_{\sigma(1)}^{(1)}, \bar z_{\tau(1)}^{(1)}) \cdots \rho(\bar
v_{\sigma(q)}^{(1)}, \bar z_{\tau(q)}^{(1)}).
$$
Write simply $\bar v_{\sigma(i)}^{(1)}=a_i, \ \bar z_i^{(1)}=b_i, \
i=1,\ldots, q$. Then
$$
\gamma_\sigma = \sum_{\tau\in S_{q}} (\mbox{sgn}\, \tau) \rho( a_1,
b_{\tau(1)}) \cdots \rho(a_{q}, b_{\tau(q)}) =\mbox{det} \left(
                \begin{array}{ccc}
                  \rho(a_1,b_1) & \cdots & \rho(a_1,b_{q}) \\
                  \vdots &  & \vdots \\
                  \rho(a_{q},b_1)& \cdots &  \rho(a_{q},b_{q}) \\
                \end{array}
              \right)
$$
\medskip
$$
= (\mbox{sgn}\, \sigma) \mbox{det}  \left(
                \begin{array}{ccc}
                  \rho(a_{\sigma^{-1}(1)},b_1) & \cdots & \rho(a_{\sigma^{-1}(1)},b_{q}) \\
                  \vdots &  & \vdots \\
                  \rho(a_{\sigma^{-1}(q)},b_1)& \cdots &  \rho(a_{\sigma^{-1}({q})},b_{q}) \\
                \end{array}
              \right) = (\mbox{sgn}\, \sigma)  \mbox{det}
              \bar\rho_1.
              $$
Hence
$$
\gamma = \frac{1}{q!}\sum_{\sigma\in S_{q}} (\mbox{sgn}\, \sigma)
\gamma_\sigma = \mbox{det} \bar\rho_1
$$
and $\varphi(g_1)= \mbox{det} \bar\rho_1\cdot \varphi(f)$. Thus,
since $\rho$ is a non degenerate form, $g_1$ does not vanish in $L$.
This completes the proof in case $k=1$.

If $k>1$, by the inductive hypothesis there exists a multilinear
polynomial
$$
g_{k-1}(x_1, \ldots, x_{q}, v_1^{(1)}, z_1^{(1)}, \ldots,
v_{q}^{(1)}, z_{q}^{(1)}, \ldots, v_1^{(k-1)},z_1^{(k-1)}, \ldots,
v_{q}^{(k-1)}, z_{q}^{(k-1)}, Y)
$$
satisfying the conclusion of the theorem. We now write
$$
g_{k-1}=g_{k-1}(x_1,\ldots, \ldots, x_{q}, Y'),
$$
where $Y'=Y \cup\{v_1^{(1)}, z_1^{(1)}, \ldots, v_{q}^{(1)},
z_{q}^{(1)}, \ldots, v_1^{(k-1)},z_1^{(k-1)}, \ldots, v_{q}^{(k-1)},
z_{q}^{(k-1)}\}$ and we apply Lemma \ref{n5} and the previous
arguments to $g_{k-1}$ . In this way we can construct the polynomial
$g_k$ and, for any evaluation $\varphi$, we have
$$
\varphi(g_k)= \mbox{det} \bar\rho_k\cdot \varphi(g_{k-1})=
\mbox{det} \bar\rho_1 \cdots \mbox{det} \bar\rho_k\cdot \varphi(f),
$$
where
$$
\bar \rho_s= \left(
  \begin{array}{ccc}
    \rho(\bar v_1^{(s)}, \bar z_1^{(s)}) & \cdots & \rho(\bar v_1^{(s)}, \bar z_{q}^{(s)}) \\
    \vdots &  & \vdots \\
    \rho(\bar v_{q}^{(s)}, \bar z_1^{(s)}) & \cdots & \rho(\bar v_{q}^{(s)}, \bar z_{q}^{(s)}) \\
  \end{array}
\right),
$$
for all $1\le s\le k$.
 This completes the proof of the proposition.
\hfill $\Box$\bigskip

\section{PI-exponents of simple color Lie algebras}

For computing PI-exponents of simple color Lie algebras we need to
get a reasonable lower bound of codimension growth.

\begin{proposition}\label{p2}
Let $L$ be as in the previous section. Then for all $n\ge 1$, there
exist constants $C>0$ and $t$ such that
\begin{equation}\label{b3}
Cn^{t}q^{n} \le c_n(L),
\end{equation}
where $q=\dim L$.
\end{proposition}

{\em Proof}. The main tool for proving the inequality (\ref{b3}) is
the representation theory of symmetric groups. We refer reader to
\cite{JK} for details of this theory.

Recall that $P_m$ is a subspace of $F\{X\}$ consisting of all
multilinear polynomials in $x_1,\ldots, x_m$ and $\Id(L)$ is the
ideal of all multilinear identities of $L$ of degree $m$. One can
define the $S_m$-action on $P_m$ b setting
$$
\sigma f(x_1,\ldots,x,m)=f(x_{\sigma(1},\ldots, x_{\sigma(m)}.
$$
Then $P_m$ becomes an $F[S_m]$-module and $P_m(L)=P_m/P_m\cap
\Id(L)$ is its submodule. By Mashke's Theorem $P_m(L)$ is the direct
sum of irreducible components and for proving the inequality
(\ref{b3}) it is sufficient to find at least one irreducible
component with the dimension greater or equal to $Cn^t q^n$.
Slightly modifying this approarch we first consider $P_{n+N}(L)$
where $n=(2k+1)q+N$ and $k,N$ are as in Proposition \ref{p1}.

Recall that there exists 1-1 correspondence between isomorphism
classes of irreducible $S_n$-representations and partitions of $n$
(or Young diagrams with $n$ boxes). A partition $\lambda\vdash n$ is
an ordered set of integers $\lambda=(\lambda_1,\ldots, \lambda_t)$
satisfying $\lambda_1\ge\ldots\ge\lambda_t>0$ and
$\lambda_1+\cdots+\lambda_t=n$. The corresponding Young diagram
$D_\lambda$ is a tableau with $n$ boxes. The first row of
$D_\lambda$ contains $\lambda_1$ boxes, the second row contains
$\lambda_2$ boxes, and so on. Young tableau $T_\lambda$ is the
diagram $D_\lambda$ filled up by integers $1,\ldots,n$.

Given a Young tableau $T_\lambda$ of shape $\lambda\vdash n$, let
$R_{T_\lambda}$ and $C_{T_\lambda}$ denote the subgroups of $S_n$
stabilizing the rows and the columns of $T_\lambda$, respectively.
If we set
$$
\bar R_{T_\lambda}= \sum_{\sigma\in R_{T_\lambda}}\sigma\quad
\mbox{and}\quad \bar C_{T_\lambda}= \sum_{\tau\in
C_{T_\lambda}}(\mbox{sgn} \tau)\tau.
$$
then the element $e_{T_\lambda}= \bar R_{T_\lambda}\bar
C_{T_\lambda}$ is an essential idempotent of the group algebra
$FS_n$ (i.e. $e_{T_\lambda}^2=\gamma e_{T_\lambda}$ for some
$0\ne\gamma\in F$) and $F[S_n]e_{T_\lambda}$ is an irreducible left
$F[S_n]$-module associated to $\lambda$.

By Proposition \ref{p1}, for any fixed $k\ge 1$ there exists a
multilinear  polynomial
$$
g_k=g_k(x_1^{(1)},\ldots, x_{q}^{(1)}, \ldots, x_1^{(2k+1)},\ldots,
x_{q}^{(2k+1)}, y_1, \ldots, y_N)
$$
such that $g_k$ is alternating on each set of indeterminates
$\{x_1^{(i)}, \ldots, x_{q}^{(i)}\}$, \ $1\le i\le 2k+1$, and $g_k$
is not a polynomial identity of $L$. Rename the variables and write
$$
g_k=h(x_1,\ldots, x_{q(2k+1)}, Y),
$$ where
$Y=\{y_1,\ldots,y_N\}$.

Since $h\not\in \Id(L)$, there exists a partition
$\lambda=(\lambda_1,\ldots, \lambda_m)\vdash n$ and a Young tableau
$T_\lambda$ such that $F[S_n] e_{T_\lambda}h \not\subseteq \Id(L)$.
Our next goal is to show that $\lambda =((2k+1)^{q})$ is a rectangle
of width $2k+1$ and height $q$.

If $\lambda_1\ge 2k+2$, then $e_{T_\lambda}h$ is a polynomial
symmetric in at least $2k+2$ variables among $x_1,\ldots, x_{n}$.
But for any $\sigma \in \bar R_{T_\lambda}$  these variables in
$\sigma \bar C_{T_\lambda}$ are divided into $2k+1$ disjoint
alternating subsets. It follows that $\sigma \bar C_{T_\lambda}h$ is
alternating and symmetric in at least two variables and so,
$e_{T_\lambda}h=0$ is the zero polynomial, a contradiction. Thus
$\lambda_1\le 2k+1$.

Suppose now that $m\ge q+1$. Since the first column of $T_\lambda$
is of height at least $q+1$, the polynomial $\bar C_{T_\lambda}h$ is
alternating in at least  $q+1$ variables among $x_1,\ldots, x_{n}$.
Since $\dim L =q$ we get that for any $\sigma$,\ $\sigma\bar
C_{T_\lambda}h\equiv 0$ on $L$ and so, also $e_{T_\lambda}h =\bar
R_{T_\lambda}\bar C_{T_\lambda}h \equiv 0$ on $L$, a contradiction.

We have proved that $F[S_n] e_{T_\lambda}h \not\subseteq \Id(L)$,
for some Young tableau $T_\lambda$ of shape $\lambda =((2k+1)^{q})$.
From the Hook formula for dimensions of irreducible representations
of $S_n$ (see \cite{JK}) and Stirling formula for factorials it
follows that
$$
\dim F[S_n]  e_{T_\lambda}h \ge \frac{q!}{(2\pi n)^q} q^n.
$$

It easily follows from the simplicity of tensor factor $B$ in
$L=F[G]\otimes B$ that
\begin{equation}\label{b4}
c_{n'}(L)\ge c_n(L) \quad \mbox{as soon as}\quad n'>n.
\end{equation}
Hence
$$
c_m(L)\ge \frac{C'}{(m-N)^q}q^{m-N}\ge \frac{C''}{m^q}q^{m}
$$
for some constants $C',C''$ for any $m=q(2k+1)+N$, $k=1,2,\ldots~$.
Finally, applying again the inequality (\ref{b4}) we get (\ref{b3})
for all $n$. \hfill $\Box$ \bigskip

Combining the inequality (\ref{ne1}) and Proposition \ref{p2} we
immediately obtain the main result of the paper.

\begin{theorem} \label{t1}
Let $F$ be an algebraically closed field of characteristic zero and
let $L=F[G]\otimes B$ be a finite dimensional color Lie superalgebra
over $F$ where $G=\<a\>_2\times\<b_2\>\simeq {\mathbb Z}_2\oplus
{\mathbb Z}_2$ with the skew-symmetric bicharacter $\beta$ defined
by $\beta(a)=\beta(b)=1, \beta(a,b)=-1$ and $B$ is a finite
dimensional simple Lie algebra with the trivial $G$-grading. Then
the PI-exponent of $L$ exists and $\mbox{exp}(L)=\dim L$.
\end{theorem}

\section{Graded identities of simple color Lie algebras}

In conclusion we discuss codimensions behavior of algebras defined
distinct bicharacters and asymptotics of graded codimensions. We
begin by an easy remark.

\begin{remark}\label{r2} If $L=F[G]\otimes B$ and $G\simeq{\mathbb Z}_2\oplus
{\mathbb Z}_2$ with the trivial bicharacter $\beta$, that is
$\beta\equiv 1$, then $\mbox{PI-exp}(L)=d=\frac{1}{4}\dim L$, where
$d=\dim B$.
\end{remark}
{\em Proof}. Since $F[G]$ is a commutative ring in this case, L is
an ordinary Lie algebra with the same identities as $B$. In
particular, $\mbox{PI-exp}(L)=\mbox{PI-exp}(B)=\dim B$ (see
\cite{GRZ} or \cite{Z}). \hfill $\Box$\bigskip

Remark \ref{r2} shows that ordinary codimensions behavior strongly
depends on bicharacter $\beta$. On the other hand one can consider
graded identities of $L$ since $L$ is a $G$-graded algebra.

Recall that if we define $G$-grading on an infinite generating set
$Y$, i.e. split $Y$ to a disjoint union $Y=\cup_{g\in G} Y^g$, $\deg
y^g=g\quad\forall~ y^g\in Y^g$, then $F\{Y\}$ can be endowed by the
induced grading if we set $\deg (y_{i_1}\cdots y_{i_m}) = \deg
y_{i_1}\cdots \deg y_{i_m}$ for any arrangement of brackets. The
polynomial $f(y^{g_1}_1,\ldots, y^{g_n}_n)$ is called a graded
identity of $L$ if $f(u_1,\ldots, u_n)=0$, as soon as $\deg u_1=\deg
y^{g_1}_1=g_1, \ldots, \deg u_n=\deg y^{g_n}_n=g_n$.

Since $F\{Y\}$ is a graded algebra, the subspace of multilinear
polynomials $P_n$ should be replaced by a graded subspace
$$
\bigoplus_{k_1+\cdots +k_4=n} P_{k_1,k_2,k_3,k_4}
$$
where $P_{k_1,k_2,k_3,k_4}$ is a subspace of  multilinear
polynomials $f$ on
$$
y_1^{g_1},\ldots,y_{k_1}^{g_1},\ldots,
y_1^{g_{k_4}},\ldots,y_{k_4}^{g_{k_4}}.
$$
Graded codimensions are defined as
$$
c_n^{gr}(L)=\sum_{{k_1\ge 0,k_2\ge 0,k_3\ge 0,k_4\ge 0\atop
k_1+k_2+k_3+k_4=n}} {n\choose k_1,k_2,k_3,k_4}
\dim\frac{P_{k_1,k_2,k_3,k_4}}{P_{k_1,k_2,k_3,k_4}\cap \Id(L)}
$$
(see \cite{GZBook} for details). For our class of algebras graded
codimensions behavior does not depend on bicharacter $\beta$
defining color on $L=F[G]\otimes B$. In the proof of the next result
we shall use the following easy observation.

\begin{remark}\label{r3} If $Lie(X)$ is a free Lie algebra on the
countable set of generators and $B$ is an arbibrary Lie algebras
then
$$
\frac{P_n}{P_n\cap \Id(B)}=\frac{V_n}{V_n\cap Id^{Lie}(B)}
$$
where $V_n$ is a subspace of $Lie(X)$ of all multilinear polynomials
in variables $x_1,\ldots, x_n$ and $Id^{Lie}(B)$ is the ideal of Lie
identities of $B$ in $Lie(X)$.
\end{remark}
\hfill $\Box$ \bigskip

We need this remark since all previous results concerning
codimension growth of Lie algebras were proved for Lie codimensions.

\begin{theorem} \label{t2}
Let $F$ be an algebraically closed field of characteristic zero and
let $L=F[G]\otimes B$ be a finite dimensional color Lie superalgebra
over $F$ where $G=\<a\>_2\times\<b_2\>\simeq {\mathbb Z}_2\oplus
{\mathbb Z}_2$ with any skey-symmetric bicharacter $\beta$ and $B$
is a finite dimensional simple Lie algebra with the trivial
$G$-grading. Then graded PI-exponent of $L$
$$
\mbox{PI-exp}^{gr}(L)=\lim_{n\to\infty}\sqrt[n]{c_n^{gr}(L)}
$$
exists and is equal to $|G|\dim L=4\dim L$.
\end{theorem}

{\em Proof}. First we prove that
\begin{equation}\label{gr1}
\dim\frac{P_{k_1,k_2,k_3,k_4}}{P_{k_1,k_2,k_3,k_4}\cap \Id(L)}=
c_n(B).
\end{equation}

Note that any multilinear Lie polynomial in $x_1,\ldots,x_n$ can be
written as a linear combination of left-normed monomials
$$
m_\sigma=[x_1,x_{\sigma(2)},\ldots,x_{\sigma(n)}],
$$
where $\sigma$ is a permutation of $2,\ldots,n$. If $w_1,\ldots,w_n$
are homogeneous elements of the color Lie superalgebra
$L=F[G]\otimes B$ then any multilinear polynomial expression on
$w_1,\ldots,w_n$ is also a linear combination of left-normed
products $[w_1,w_{\sigma(2)},\ldots,w_{\sigma(n)}]$. Since we are
interested in graded identities of $L$ it is sufficient to consider
only left-normed monomials and their linear combinations.

Let $(g_1,\ldots,g_n)$ be a n-tuple of elements of $G$ such that
$$
(g_1,\ldots,g_n)=(\underbrace{e,\ldots,e}_{k_1},
\underbrace{a,\ldots,a}_{k_2},\underbrace{b,\ldots,b}_{k_3},
\underbrace{ab,\ldots,ab}_{k_4}).
$$
Then
$$
m_\sigma(g_1\otimes x_1,\ldots,g_n\otimes x_n)= g_1\cdots g_n\otimes
\lambda_\sigma m_\sigma(x_1,\ldots,x_n)
$$
in $F[G]\otimes F\{X\}$ where $\lambda_\sigma=\pm 1$ depends only on
$\sigma$ for given $g_1,\ldots, g_n$.

Given a multilinear polynomial
$$
f=f(x_1,\ldots,x_n)=\sum_{\sigma\in S_{n-1}} \alpha_\sigma m_\sigma
$$
of $F\{X\}$ we denote by $\widetilde f$ the element
$$
\widetilde f=\widetilde f(x_1,\ldots,x_n)=\sum_{\sigma\in S_{n-1}}
\lambda_\sigma\alpha_\sigma m_\sigma.
$$
Then for any $w^1_1,\ldots,w^1_{k_1},\ldots,
w^4_1,\ldots,w^4_{k_4}\in B$ we have
$$
f(e\otimes w^1_1,\ldots, e\otimes w^1_{k_1},\ldots, ab\otimes
w^4_1,\ldots, ab\otimes w^4_{k_4})
$$
$$
=a^{k_2}b^{k_3}(ab)^{k_4}\otimes \widetilde
f(w^1_1,\ldots,w^1_{k_1},\ldots, w^4_1,\ldots,w^4_{k_4}).
$$

In particular, $f$ is a graded identity of $L$, $f\in
P_{k_1,k_2,k_3,k_4}\cap \Id(L)$ if and only if $f$ is an identity of
the Lie algebra $B$. Now, if $c_n(B)=N$ and $m_{\sigma_1}, \ldots,
m_{\sigma_N}$ is a basis of $V_n$ in $Lie(X)$ modulo $Id^{Lie}(B)$
then also $m_{\sigma_1}, \ldots, m_{\sigma_N}$ is a basis of
$P_{k_1,\ldots, k_4}$ modulo $\Id(L)$ in $F\{X\}$ and we have proved
the relation (\ref{gr1}). Hence
$$
c_n^{gr}(L)=\sum_{{k_1\ge 0,k_2\ge 0,k_3\ge 0,k_4\ge 0\atop
k_1+k_2+k_3+k_4=n}} {n\choose k_1,k_2,k_3,k_4}
\dim\frac{P_{k_1,k_2,k_3,k_4}}{P_{k_1,k_2,k_3,k_4}\cap \Id(L)}
$$
$$
=c_n(B)\sum_{{k_1\ge 0,k_2\ge 0,k_3\ge 0,k_4\ge 0\atop
k_1+k_2+k_3+k_4=n}} {n\choose k_1,k_2,k_3,k_4}=4^n c_n(B)
$$
and we have completed the proof since
$\lim_{n\to\infty}\sqrt[n]{c_n(B)}=\dim B$ by \cite{GRZ}. \hfill
$\Box$\bigskip

The first and the second author were supported by the Slovenian
Research Agency grants P1-0292-0101-04 and J1-9643-0101. The third
author was partially supported by RFBR grant No
 09-01-00303  and SSC-1983.2008.1.


\end{document}